\documentclass{amsart}
\usepackage{graphicx,amssymb}
\newtheorem{theorem}{Theorem}
\newtheorem{proposition}[theorem]{Proposition}
\newtheorem{lemma}[theorem]{Lemma}
\newtheorem{corollary}[theorem]{Corollary}
\newtheorem{remark}[theorem]{Remark}
\newtheorem{example}[theorem]{Example}

\newcommand{\bth}{\begin{theorem}}
\newcommand{\bpr}{\begin{proposition}}
\newcommand{\epr}{\end{proposition}}
\newcommand{\bco}{\begin{corollary}}
\newcommand{\eco}{\end{corollary}}
\newcommand{\ble}{\begin{lemma}}
\newcommand{\ele}{\end{lemma}}
\newcommand{\bre}{\begin{remark}\rm}
\newcommand{\ere}{\end{remark}}
\newcommand{\bex}{\begin{example}\rm}
\newcommand{\eex}{\end{example}}

\def\la#1{\hbox to #1pc{\leftarrowfill}}
\def\ra#1{\hbox to #1pc{\rightarrowfill}}
\def\fract#1#2{\raise4pt\hbox{$ #1 \atop #2 $}}
\def\decdnar#1{\phantom{\hbox{$\scriptstyle{#1}$}}
\left\downarrow\vbox{\vskip10pt\hbox{$\scriptstyle{#1}$}}\right.}

\def\lrar{{\ra 2}}
\def\tensor{\otimes}
\def\sp#1{\hbox{SP}^{#1}}
\def\bsp#1{\overline{\hbox{SP}}^{#1}}
\def\spy{\hbox{SP}^{\infty}}
\def\bspy{\overline{\hbox{SP}}^{\infty}}

\def\bbz{{\mathbb Z}}

\def\bbp{{\mathbb P}}
\def\bbr{{\mathbb R}}

\def\bbc{{\mathbb C}}
\def\bbr{{\mathbb R}}
\def\bbq{{\mathbb Q}}
\def\bbn{{\mathbb N}}

\begin{document}

\title{Symmetric Products of Two Dimensional Complexes}

\author{Sadok Kallel}

\address{Laboratoire Painlev\'e, Universit\'e de Lille I, Villeneuve d'Ascq,
France}
\email{sadok.kallel@math.univ-lille1.fr}

\author{Paolo Salvatore}

\address{Dipartimento di matematica, Universit\`a di Roma ``Tor Vergata'',
Italy}
\email{salvator@mat.uniroma2.it}


\begin{abstract} We exhibit a multiplicative and minimal cellular
complex which allows explicit and complete (co)homological
calculations for the symmetric products of a finite two dimensional CW
complex. By considering cohomology, we observe that a classical
theorem of Clifford on the dimension of various linear series on a
projective curve has a purely topological statement. We give a ``real"
analog of this theorem for unoriented topological surfaces.
\end{abstract}

\maketitle

\centerline{To Sam Gitler on the occasion of his 70th birthday}


\section{One dimensional Complexes; Constructions}

The $n$-th symmetric product $\sp{n}X$ of a space $X$, which is
assumed to be connected and based, is the quotient of $X^n$ by the
permutation action of the symmetric group. We write an element of
$\sp{n}X$ as an unordered $n$-tuple of points $<x_1,\ldots, x_n>$.

When $X=S^1$ there are several ways to see that $\sp{n}(S^1)$ is up to
homotopy the circle again. An easy such way is to replace up to
homotopy $S^1$ with the punctured complex plane $\bbc^*$.  Then the
map
$$ \sp{n}(\bbc^*)\lrar {\mathcal H},\ \
< z_1,\ldots,z_n> \mapsto\prod_{1\leq i\leq n}(z-z_i) 
$$
is a homeomorphism, where on the right we have identified the space of
monic polynomials with roots avoiding the origin with $\mathcal H$,
the complement of a hyperplane in $\bbc^n$. But then evidently
${\mathcal H} \cong\bbc^{n-1}\times\bbc^*\simeq S^1$ and hence the
claim.  It is amusing to see for instance that $\sp{2}(S^1)$ is the
Moebius band. 
In general one has the following
complete description due to Morton.

\bth\label{sym1}\cite{morton}
For all $n\geq 1$, the multiplication map $$m : \sp{n}S^1\lrar S^1\ \
,\ \ <x_1,\ldots, x_n>\mapsto x_1\cdots x_n$$ has the
structure of an $(n-1)$-dimensional disc bundle over $S^1$; trivial if
$n$ is odd and non-oriented if $n$ is even.
\end{theorem}

As is clear from theorem \ref{sym1}, a chain complex for the symmetric
product of a circle has a small retract.  This suggests the following
construction.

Let $X$ be any reduced CW complex of finite type, and let $\bigvee^k
 S^1\hookrightarrow X$ be the one-skeleton inclusion. We choose the
 basepoint $*\in X$ to be the unique 0-cell and we assume that this
 basepoint corresponds to the identity element $1$ in each circle leaf
 viewed as a circle group.

Define the identification space
$$\bsp{n}X = \sp{n}X/\sim$$
where $\sim$ identifies $<x,y, z_1,\ldots, z_{n-2}>$ to
$<*,xy, z_1,\ldots, z_{n-2}>$
whenever $x, y$ are in the same leaf $S^1$. We denote by $q$ the
quotient map.  The usefulness of this construction is contained in the
following crucial lemma

\ble The projection $q:\sp{n}X\to\bsp{n}X$ \label{1} is a homotopy
 equivalence.  \ele

\begin{proof}
We start by the case $X=X^{(1)}$, i.e. when $X$ is a bouquet of $k$
circles.  The space $\sp{n}(X)$ is the colimit of the diagram
$(i_1,\dots,i_k) \mapsto \sp{i_1}(S^1) \times \dots \times
\sp{i_k}(S^1)$, indexed over the full sub-poset $I \subset \bbn^k$
containing $k$-tuples with $i_1+\dots+i_k \leq n$.  The morphisms in
the diagram are induced by the obvious inclusions.  The space
$\bsp{n}(X)$ is the colimit of a similar diagram sending
$(i_1,\dots,i_k)$ to $(S^1)^{e(i_1)} \times \dots \times
(S^1)^{e(i_k)} $, where $e(x)=1$ if $x>0$ and $e(x)=0$ if $x=0$.  The
circle multiplication induces on colimits the projection $q$.  It is
not difficult to see (compare \cite{DS}) that each colimit is homotopy
equivalent to the associated homotopy colimit. The circle
multiplication being an equivalence on each factor by Morton's result,
the lemma follows in this case by homotopy invariance.

In general we have a filtration
$E_0 \subset  \dots \subset  E_n =\sp{n}(X)$
where $E_i$ is the space of those
unordered $n$-tuples containing at most $i$ points outside the 1-skeleton.
Each subspace $E_i$ is closed in $E_{i+1}$, it is the strong deformation retract
of an open neighbourhood $U_i$ in $E_{i+1}$, and all this passes to the
quotient. For example choose $U_i$ to contain at most $i$ points outside
an open neighbourhood $U$ of $X^{(1)}$ in $X$, with $X^{(1)}$ deformation
retract of $U$.
The projection between differences $E_{i}-E_{i-1} \to q(E_{i})-q(E_{i-1})$
can be identified to the projection
$$\sp{i}(X-X^{(1)}) \times \sp{n-i}(X^{(1)})
\to \sp{i}(X-X^{(1)}) \times \bsp{n-i}(X^{(1)})$$
 and is a homotopy equivalence by our earlier remark. We can now
use the gluing lemma and
induction up the filtration to show that $E_i\simeq q(E_i)$,
$1\leq i\leq n$.
\end{proof}

As it turns out, the construction $\bsp{}$ is much easier to study.
For example it is straightforward to see that

\ble $\bsp{n}(\bigvee^kS^1 )$ is homeomorphic to the $n$-skeleton of
$(S^1)^k$.\ele

As an immediate corollary we obtain the following main result of
\cite{ong} which initially made a lengthy use of the theory of
hyperplane arrangements.

\bco\label{circles}
For $n\geq k\geq 1$
 there is a homotopy equivalence
$\sp{n}(\bigvee^kS^1 ) \simeq (S^1)^k$. For
$n < k$  $\sp{n}(\bigvee^kS^1 )$ is homotopic to the union of
$k\choose n$ $n$-dimensional subtori in $(S^1)^k$.  \eco


\section{Two Dimensional Complexes and Minimal Cell Decompositions}
\label{2complexes}

Let $X$ be a connected (based) cellular complex
\begin{equation}\label{twocell}
X = (\bigvee^k S^1)\cup (D^2_1\cup\cdots\cup D^2_r)
\end{equation}
obtained by attaching $r$ 2-cells to a bouquet of $k$ circles. A chain
complex for $X$ has 1-dimensional generators $e_1,\ldots, e_k$ and two
dimensional generators $D_1,\ldots, D_r$.  In this section we show
that these classes generate multiplicatively (in a sense we define
shortly) a minimal cell complex for $\sp{n}X$.

The main observation in dimension two is the following standard
result.

\ble\label{sym2} Let $\dot{D}^n$ be the open unit disc in $\bbr^n$,
 and $D^n$ its closure.  Then there is a homeomorphism of pairs
 $\phi_n:(\sp{n}D^2,\sp{n}D^2 -\sp{n} \dot{D}^2 )
 \cong(D^{2n},S^{2n-1})$.  \ele

\begin{proof}
For a positive integer $d$ we define a self homeomorphism
$r_d: \bbc \to \bbc$ by
 $r_d(z)=z|z|^{1/d-1}$ for $z \neq 0$ and
 $r_d(0)=0$. In particular if $z \geq 0$ then $r_d(z)=z^{1/d}$.
Let
$f_d<z_1,\ldots,z_n>$ be the elementary degree $d$ symmetric function
in $z_1,\ldots,z_n$.
By the fundamental theorem of algebra
the correspondence
$$<z_1,..., z_n>\mapsto (f_1, ..., f_n)$$
gives a homeomorphism $SP^n(\bbc) \cong \bbc^n$.
Let us write
$\partial \sp{n} D^2 =
 \sp{n}D^2 -\sp{n}\dot{D}^2$.
Via the action $t<z_1,\ldots,z_n>=<tz_1,\ldots,tz_n>$, $t \geq 0$, we may
identify $\sp{n} \bbc =\sp{n} \bbr^2$ to the $\bbr_+$-cone on
$\partial \sp{n}D^2$ and
$\sp{n}D^2$ to the $I$-cone on $\partial \sp{n}D^2$.
Of course ${D}^{2n}$ and $\bbr^{2n}$ are respectively the $I$-cone and
the $\bbr_+$-cone on $S^{2n-1}$.
It is easy to see that the homeomorphism $\psi_n:SP^n(\bbc) \cong \bbc^n$
sending $<z_1,\ldots,z_n> \mapsto (r_1 (f_1) , \ldots,r_n (f_n) )$
is $\bbr_+$-equivariant and
preserves cone rays.
On $\partial \sp{n}D^2$ we define $\phi_n=\psi_n/\|\psi_n\|$
and extend it to $\sp{n}D^2$ by cone extension.
\end{proof}

A quick corollary of lemma \ref{sym2} is that $\sp{n}S$ is an
orientable (resp. non-orientable) $2n$-dimensional manifold if $S$ is
an orientable (resp. non-orientable) topological surface
(eg. \cite{gunning}).

Write $\bsp{}X = \coprod_{n\geq 0}\bsp{n}X$ where $\bsp{0}X =
\centerdot$ is the basepoint and $\bsp{1}X = X$.  Note that
concatenation $(S^1)^s\times (S^1)^t\lrar (S^1)^{s+t}$ commutes with
multiplication in the circle group $S^1$ so that we have an induced
pairing at the level of quotient spaces
\begin{equation}\label{diagram}
\begin{matrix}
X^s\times X^t&\fract{}{\lrar}& X^{s+t}\cr
\decdnar{\bar q_s\times \bar q_t}&&\decdnar{\bar q_{s+t}}\cr
\bsp{s}X\times\bsp{t}X&\fract{*}{\lrar}&\bsp{s+t}X\cr
\end{matrix}
\end{equation}
This endows $\bsp{}X$ with the structure of a commutative monoid.
A cellular decomposition of $\bsp{}X$ is called \emph{multiplicative} if
it is compatible with this monoid structure; i.e. if the multiplication
$*$ is cellular. Such a decomposition specializes to
a decomposition for $X$ such that the projections
\begin{equation}\label{firstquotient}
\bar q : X^n\lrar\bsp{n}X\ \ ,\ n\geq 1
\end{equation}
are cellular. Given a multiplicative decomposition for $\bsp{}X$;
$c_1\in C_*(\bsp{s}X)$ and $c_2\in C_*(\bsp{t}X)$, we denote
by \emph{the product}
$c_1*c_2$ the image of $c_1\tensor c_2$ under the map
$C_*(\bsp{s}X)\tensor C_*(\bsp{t}X)\lrar C_*(\bsp{s+t}X)$.

\bth\label{main1} Write $X = \bigvee^kS^1\cup (D^2_1\cup\cdots\cup D^2_r)$. Then\\
(a) A chain complex for $\bsp{}X$ is multiplicatively generated under $*$
by  a zero dimensional class $v_0$, degree one classes $e_1,\ldots, e_k$
and degree $2s$ classes
$\sp{s}D_i$ $1\leq i\leq r$, $1\leq s$, under the relations
\begin{eqnarray*}
e_i*e_j = -e_j*e_i \, (i \neq j)\ \ ,\ \ e_i*e_i = 0 \ \ \\
\sp{s}D_i*\sp{t}D_i = {s+t\choose t}\sp{s+t}D_i \,
\end{eqnarray*}
The boundaries are such that:
\begin{eqnarray*}
&&\partial e_j = 0\ \ ,\ \
\partial\sp{s}D_i = (\partial D_i)*\sp{s-1}D_i\\
&&\partial\ \ \hbox{is a derivation}
\end{eqnarray*}
where $\partial D_i$ is the boundary of $D_i$ in $C_1(X)$.\\
(b) A cellular chain complex for $\bsp{n}(X)$ consists of the subcomplex
generated by
$$v_0^k*e_{i_1}*\cdots *e_{i_t}*\sp{s_1}(D_{j_1})*\cdots *\sp{s_l}(D_{j_d})$$
with $k+t+s_1+\cdots +s_l= n$.\\
(c) $H_*(\sp{n}X;\bbz )$ and
$H_*(X;\bbz )$ have the same prime torsion.\\
(d) In particular $H_*(X)$ is torsion free if and only if
$H_*(\sp{n}X)$ is torsion free, in which case there is a commutative
diagram
$$\begin{matrix} H_*(\bsp{n}X)&\hookrightarrow&C_*(\bsp{n}X)\cr
\decdnar{}&&\decdnar{}\cr H_*(\bsp{}X)&\hookrightarrow&C_*(\bsp{}X)\cr
\end{matrix}$$
with all maps being monomorphisms.  The bottom map in the diagram
above is a ring map with respect to the symmetric product *.
\end{theorem}

\begin{proof} 
Start with $X$ having one skeleton $\bigvee^kS^1\hookrightarrow X$
with zero-cell the wedge point, and two cells attaching to the
bouquet.  Write the zero cell $\centerdot$, the one dimensional cells
$E_1,\ldots, E_k$ and the two dimensional cells $D_1^2,\ldots, D_r^2$.
These cells are closed with disjoint interiors homeomorphic to discs.

Suppose $C_1$ and $C_2$ are cells of such a CW decomposition of $X$,
$q : X\times X\lrar\sp{2}X$ the quotient map.  The image $q(C_1\times
C_2)$ is generally not a cell in $\sp{2}X$. We exploit the fact that
this problem doesn't occur when we map into $\bsp{2}X$.

Write $\bar q: X^n \lrar \bsp{n}X$
the quotient map, and
denote by the product $C_1*\ldots *C_n$
the image of $C_1\times \ldots \times  C_n$ under $\bar q$.
We also write
$\bar q(\centerdot)=v_0$, and identify 1-cells and 2-cells with their
image via $\bar q$.
We use the following geometric properties:           \\
(i) The $*$-product of cells of $X$ whose interiors are disjoint 
is again a cell.
This is because $\bar{q}$ is injective on the product of the interiors.\\
(ii) $E_i*E_i = E_i*v_0$ in $\bsp{2}X$.\\
(iii) The $n$-fold product of $D^2_s$ via $*$
 is a cell that is covered $n!$ times by
$(D^2_s)^n$ according to lemma \ref{sym2}. We will
denote by $\sp{n}(D^2_s)$ this symmetric product cell.
\

This yields a CW decomposition of
$\bsp{n}X$ for all $n\geq 1$ such that $\bar q
:X^n\rightarrow\bsp{n}X$ is cellular. In other words,
a multiplicative CW decomposition
for $\bsp{}X:=\coprod_{n\geq 0}\bsp{n}X$ is obtained from $v_0$,
the $E_i$'s, the $D_j^2$'s and all possible products under $*$ among
these with the restriction that there are only a
finite number of non-zero cells in the product. If we agree to identify
cells when 
 they differ by multiples of $v_0$, then the CW structure on
$\bsp{n}X$ is obtained by taking all $m$-fold products with $m\leq n$.

We can now pass to the chain complex level. To $E_i$ corresponds the
algebraic generator $e_i\in C_1(\bsp{n}X)$, to $D^2_j$ corresponds
$D_j$ of degree two with well defined boundary.  This determines
$C_*(X)$. We extend this complex multiplicatively to all of
$C_*(\bsp{}X)$ using the diagram obtained from (\ref{diagram})
\begin{equation}\label{diagram2}
\begin{matrix}
C_i(X^s)\otimes C_j(X^t)&\fract{}{\lrar}& C_{i+j}(X^{s+t})\cr
\decdnar{\bar q_*\otimes \bar q_*}&&\decdnar{\bar q_*}\cr
C_i(\bsp{s}X)\otimes C_j(\bsp{t}X)&\fract{*}{\lrar}&C_{i+j}(\bsp{s+t}X)\cr
\end{matrix}
\end{equation}

First of all, since $E_i*E_i = E_i$, we set $e_i*e_i = 0$.  For chains
$c_1$ and $c_2$ supported by geometric cells $C_1$ and $C_2$ with
disjoint interiors, the diagram implies that $c_1*c_2$ is the cell
supported by $C_1*C_2$.  More interestingly note that $\bar
q_*(D_i^{\tensor s}) = s!\sp{s}(D_i)$ so that tracing through
(\ref{diagram2}) we must have that $\sp{s}(D_i)*\sp{t}(D_i) =
{s+t\choose t}\sp{s+t}D_i$.  Finally by commutativity we have that
$c_i*c_j = (-1)^{|c_i||c_j|}c_j*c_i$.

We next analyze the boundaries.  It is clear that $\partial e_i= 0$
while $\partial D_j$ is determined by the attaching maps of $X$. We
need understand the boundary on the ``symmetric product" cell
$\sp{s}(D)$. This is described geometrically as the image under the
symmetric quotient $X^s\lrar\sp{s}X$ of
\begin{equation}\label{boundary}
\partial (D^2_i)^s =
\bigcup_{j=1} (D^2_i)^{j-1}\times\partial D^2_i\times (D^2_i)^{s-j}
\end{equation}
and each term in the union maps to $(\partial D^2_i)*\sp{s-1}(D^2_i)$
(in an orientation preserving manner).  But the degree of the
projection of the right hand side of (\ref{boundary}) into $(\partial
D^2_i)*\sp{s-1}(D^2_i)$ is $s\cdot (s-1)! = s!$ (taking into account
$s$-terms in the union and then the degree of the projection
$(D^2_i)^{s-1}\lrar\sp{s-1}(D^2_i)$). On the other hand $s!$ is
precisely the degree of the projection $\partial
(D^2_i)^s\lrar\partial \sp{s}(D^2_i)$ so that in the chain complex
we must have $\partial\sp{s}(D_i) = (\partial D_i)*\sp{s-1}(D_i)$. 

We observe that $\partial$ is a derivation because the cellular
decomposition is multiplicative. The remaining claims follow by
construction.

Parts (c) and (d) follow immediately from parts (a) and (b).
\end{proof}

\bre\label{bidegree} We have embeddings
$\bsp{s-1}X\hookrightarrow\bsp{s}X$ (adjunction of basepoint) and
$C_*(\bsp{s-1}X)\hookrightarrow C_*(\bsp{s}X)$ (multiplication by
$v_0$).  We can then assign a bidegree to cells so that $c$ has
bidegree $(s,*)$ if $c\in C_*(\bsp{s}X,\bsp{s-1}X)$. We refer to $s$
as the \emph{filtration} degree.  For example $\sp{s}(D)$ has bidegree
$(s,2s)$ and the product of distinct terms $e_1*\cdots *e_s$ has
bidegree $(s,s)$.  The useful feature here is that the boundary
operator $\partial$ preserves filtration degrees so that we have a
decomposition
$$H_*(\bsp{n}X)\cong\bigoplus_{1\leq i\leq n}H_*(\bsp{i}X,\bsp{i-1}X)$$
which is a special case of a more general splitting result of Steenrod
\cite{dold}.
\ere

\bre\label{infinite} It is often convenient to consider as in the
literature the infinite symmetric product $\spy (X)$ obtained as the
direct limit of the basepoint inclusions $\sp{n}X\lrar\sp{n+1}X$. Let
$\bspy (X)$ be the induced quotient, $\bspy (X)\simeq\spy (X)$.  A
chain complex for $\bspy (X)$ is obtained from a chain complex for
$\bsp{}(X)$ by identifying cells differing by a multiple of $v_0$.
\ere

\noindent{\sc Examples}. We illustrate theorem \ref{main1} with a few
examples.
\begin{enumerate}
\item First
when $X = S^2 = *\cup D^2$ and there are no one cells. Here $\partial
D^2=0$ and the homology of $\sp{n}X = \bsp{n}X$ is generated in
dimension $2i$ by the unique cell in that dimension $\sp{i}(D)$,
$i\leq n$. This is of course in accordance with the identification
$\sp{n}S^2 = \bbc\bbp^n$.
\item If instead we write $S^2 = S^1\cup D^2_1\cup D^2_2$ with
$\partial D_1=e_1$ and $\partial D_2=-e_1$, then $H_*(\bspy (X)) =
H_*(\spy (X))$ has generator $D_1+D_2$ in dimension 2, $\sp{2}(D_1) +
D_1*D_2 + \sp{2}(D_2)$ in dimension 4 and more generally
$\sum_{s+t=n}\sp{s}(D_1)*\sp{t}(D_2)$ in dimension $2n$.
\item Write $\bbr\bbp^2 = S^1\cup_f D^2$ where the attaching map is of
degree two. A chain complex for $\bbr\bbp^2$ has generators $e,D$ with
$\partial D = 2e$. A chain complex for $\bsp{n}X$ has even generators
$\sp{i}(D)$ and odd generators $e \sp{i}(D)$ with

\begin{eqnarray*}
\partial (\sp{i}(D))
&=& 2e\sp{i-1}(D),\\
\partial (\sp{i}(D)e) &=& 2e^2\sp{i-1}D = 0
\end{eqnarray*}
One sees immediately that $H_*(\sp{n}(\bbr\bbp^2)) = H_*(\bbr\bbp^{2n})$
in accordance with the
identification $\sp{n}(\bbr\bbp^2)\cong\bbr\bbp^{2n}$ (lemma \ref{dupont}).
\item If $S$ is a closed Riemann surface, then
$H_*(\sp{n}S)$ is torsion-free for all $n\geq 1$ (classical, cf. \cite{macdo}).
\end{enumerate}

\bre Notice that the largest cells in $C_*(\bsp{n}X)$ have dimension
$2n$ and are of the form $\sp{n}(D)$ for some two cell $D\in
C_2(X)$. This implies that if $X$ is a two complex, then
$H_*(\sp{n}X;\bbz )$ is trivial for $*\geq 2n+1$. More generally if
$X$ is an $m$-dimensional complex, then $\sp{n}X$ is
$nm$-dimensional.\ere

\bco\label{abelian} $\pi_1(SP^{n}X)$ is abelian for $n>1$.  \eco

\begin{proof}
By theorem \ref{main1}, $\bsp{n}X$
and $\bspy X$ have the same two skeleton and hence the same
$\pi_1$. But $\bspy X$ is an abelian monoid and hence
has abelian fundamental group.
\end{proof}

The next two corollaries are of good use in applications and combine
remark \ref{infinite} with theorem \ref{main1} (compare \cite{jim})

\bco\label{divided} Let $X$ be a 2-complex and suppose that $\partial
D = 0$ for some two cell $D$. Then the cells $\sp{s}(D)$ for $s\geq 0$
generate a divided power algebra $\Gamma (D)$ in $H_*(\spy (X))$. \eco

\bco\label{Riemann} Let $X = (\bigvee^{2g} S^1)\cup D^2$ be a Riemann
surface of genus $g$. Then a minimal multiplicative chain complex for
$\bspy (X)$ is given by
\begin{equation}\label{thecomplex}
E(e_1,\ldots, e_{2g})\tensor\Gamma (D).
\end{equation}
At every finite stage $\bsp{n}X$ has a minimal cell complex consisting
of all chains of filtration $\leq n$ in the bigraded complex
(\ref{thecomplex}) above, where $e_i$ is of bidegree $(1,1)$ and $D$
of bidegree $(1,2)$ (remark \ref{bidegree}).  The class $\sp{n}(D)$
corresponds in homology to the orientation class of the manifold
$\sp{n}(X)$. \eco

\subsection{Comparison with methods of Dold and Milgram}
Dold \cite{dold} and then Milgram \cite{jim} gave an effective recipe
to compute the homology of symmetric products of CW complexes. The
idea is that if $X$ is any complex of the homology type of a wedge of
Moore spaces $\bigvee_{i=1}^r M_i$, then as graded abelian groups
\begin{equation}\label{formula}
\begin{split}
&H_*(\sp{n}X,\sp{n-1}X)\cong\\
&\bigoplus_{i_1+\cdots + i_r = n} H_*(\sp{i_1}M_1,\sp{i_1-1}M_1)\tensor
\cdots \tensor H_*(\sp{i_r}M_r,\sp{i_r-1}M_r)
\end{split}
\end{equation}
with the tensor product on the right corresponding to the symmetric
product pairing * in homology on the left.
For example in the case of a Riemann surface $S$ of genus $g$,
$H_*(S)\cong H_*(\bigvee^{2g}S^1\vee S^2)$ and hence according to
(\ref{formula})
we recover corollary \ref{Riemann}.

The decomposition in (\ref{formula}) however does not shed light on
neither cup products nor cohomology operations in $\sp{n}X$. We will
deal with this in the next section.


\section{Cohomology Structure}\label{cohomology}

Let $X$ be a two dimensional complex and suppose that $H_*(\sp{n}X)$
is torsion free. Then the transfer shows that
$$H^*(\sp{n}X) = \left(H^*(X)^{\tensor n}\right)^{\Sigma_n}$$
is the submodule of invariants. In particular the induced map
in cohomology $\pi^* : H^*(\sp{n}X)\lrar H^*(X^n )$ is injective.
This is the method adopted by
MacDonald in \cite{macdo} to determine the cohomology ring
of the symmetric product of an orientable surface $S$
(i.e. theorem \ref{macdocalc}).

The situation for more general $X$ is harder to track down as $\pi^*$
is not necessarily injective (eg. this is already not the case for $X
= \bbr\bbp^2$ and $n=2$).  To remedy to this problem, we need use other
arguments based on the multiplicative cell complex introduced in
section \ref{2complexes}.

Let $\Delta$ be the diagonal map and $\delta :X\lrar X\times X$ a
cellular approximation.  Write $H : \Delta\simeq\delta$ for the
homotopy. For reasons that will soon be clear, we would like to choose
$H$ so that $H_t, t\in [0,1]$ sends
 each leaf of the bouquet $\bigvee
S^1\subset X$ to its square.
  Start with a standard approximation for the
diagonal on $S^1$.  This can be done on each leaf to yield an
approximation $\delta^{\vee}$ for the diagonal $\Delta^{\vee}$ on the
bouquet. The relative cellular approximation theorem (\cite{hatcher},
theorem 4.8) states that it is possible to extend $\delta^{\vee}$ to
a cellular map $\delta$ on all of $X$.

We wish to understand the cup product of $H^*(\bsp{n}X)$ starting from

$\delta_*:C_*(X)\lrar C_*(X)\tensor C_*(X)$.  The first step is to
consider the coproduct for $C_*(X^n)$ which is obtained up to suitable
shuffle from the map $\delta_*^{\tensor n}$.  More explicitly, if
$\chi$ is the shuffle map
$$\chi : (x_1,\ldots, x_n)\times (y_1,\ldots, y_n)\mapsto
(x_1,y_1)\times (x_2,y_2)\times\ldots\times (x_n,y_n)$$
then we can write the diagonal
$(x_1,\ldots, x_n)\mapsto (x_1,\ldots, x_n; x_1,\ldots, x_n)$
as a composite of $(x_1,\ldots, x_n)\mapsto (x_1,x_1, \cdots , x_n,x_n)$
followed by $\chi^{-1}$.
A diagonal approximation for $X^n$ is then given by
$X^n\fract{\delta^n}{\lrar} (X^2)^{n} \fract{\chi^{-1}}{\lrar} X^n\times X^n$.

Suppose now that $\coprod\bsp{n}X$ is given a multiplicative cell
decomposition (as in section \ref{2complexes}) so that in particular
the quotient $\pi : X^n\lrar\bsp{n}X$ is cellular.

\ble\label{coproduct} There is a commutative diagram
$$\begin{matrix}
C_*(X)^{\tensor n}&\fract{\chi_*^{-1}\delta_*^n}{\ra 3}&C_*(X)^{\tensor n}
\tensor C_*(X)^{\tensor n}\cr
\decdnar{\pi_*}&&\decdnar{\pi_*\tensor\pi_*}\cr
C_*(\bsp{n}X)&\fract{\lambda_*}{\ra 3}&C_*(\bsp{n}X)\tensor C_*(\bsp{n}X)\cr
\end{matrix}$$
where $\lambda_*$ induces in cohomology the cup product.
\ele

\begin{proof}
The cellular approximation of the diagonal of $X^n$ induces a map

$\lambda:\sp{n}X \lrar \sp{n}X \times \sp{n}X$ homotopic to the
diagonal, but not a map $\bsp{n}X \to \bsp{n}X \times \bsp{n}X$
because the approximation on a leaf $S^1 \to S^1 \times S^1$ is not a
homomorphism.  Let us filter $\sp{n}X$ by the inverse images of the
skeleta of $\bsp{n}X$. The proof of lemma \ref{1} shows that the
projection from an inverse image of a skeleton to the skeleton is a
homotopy equivalence.  Thus the Leray spectral sequence of our
filtration has as $E_1$ term the chain complex $C_*(\bsp{n}X)$ and
collapses at the $E_2$ term. If we filter similarly the product
$\sp{n}X \times \sp{n}X$ then $\lambda$ is a filtration preserving map
inducing $\lambda_*$ on the $E_1$-term.  The commutative diagram in
the statement lives at the level of Leray $E_1$-terms, where $X^n$ and
its square are filtered by skeleta.
\end{proof}

\subsection{A calculation}
As an illustration of the method and for later use, we determine the
cohomology of $\sp{n}(S^1\cup_{m}D^2 )$, where $X=S^1\cup_{m}D^2$ is
the complex obtained by attaching $D^2$ along a degree $m$ map. The
chain complex for $X$ has generators $e, D$ with $\partial D =
me$. The cup product structure is only interesting with $\bbz_m$
coefficients.  In this case, $D$ is primitive if $m$ is odd, and
otherwise (see \cite{hatcher}, example 3.9)
\begin{equation}\label{3.1.1}
D\fract{\lambda_*}{\lrar} D\tensor 1 + ke\tensor e + 1\tensor D\ \ \ \ ,
\ \ m = 2k
\end{equation}
To determine $\lambda_*(\sp{2}D)$ we look at the coproduct upstairs
in lemma \ref{coproduct}
\begin{eqnarray*}
D\tensor D
&\mapsto&\chi^{-1}\left(( D\tensor 1\tensor 1\tensor 1 +
k e\tensor e\tensor 1\tensor 1 + 1\tensor D\tensor 1\tensor 1)\right.\\
&&\ \left. (1\tensor 1\tensor D\tensor 1 +
k 1\tensor 1\tensor e\tensor e + 1\tensor 1\tensor 1\tensor D )\right)\\
&& = D\tensor D\tensor 1\tensor 1 + kD\tensor e\tensor 1\tensor e +
D\tensor 1\tensor 1\tensor D \\
&&\ \ + ke\tensor D\tensor e\tensor 1 -k e^{\tensor 4} +
ke\tensor 1\tensor e\tensor D \\
&&\ \ +  1\tensor D\tensor D\tensor 1 + k1\tensor e\tensor D\tensor e +
1\tensor 1\tensor D\tensor D
\end{eqnarray*}

We can now apply $\pi_*$ to the left and $\pi_*\tensor\pi_*$ to the right
to obtain (after dividing by $2$)
\begin{equation}\label{3.1.2}
\sp{2}D\fract{\lambda_*}{\longmapsto}
\sp{2}D\tensor 1 + k De\tensor e + D\tensor D + ke\tensor eD
+ 1\tensor \sp{2}D
\end{equation}
Note that $H^*(\sp{n}X;\bbz_m)$ has one generator per dimension.
Denote by $b:= D^*$ the dual of $D$ and by $f = e^*$ the dual of $e$. Then
(\ref{3.1.1}) implies that $f^2 = kb$ while (\ref{3.1.2}) implies
that $b^2 = (\sp{2}D)^*$ and $f^2b = k(\sp{2}D)^* = kb^2$. Carrying this
game to the remaining classes shows that $b$ generates a truncated
polynomial algebra where $b^k$ is dual to $\sp{k}D$ and $b^{n+1}=0$.
On the other hand, $fb^k$ is dual to $e\sp{k}D$ (compare (\ref{exx})).
This yields

\ble If $m$ is even, $m=2k$, then
$H^*(\sp{n}(S^1\cup_mD);\bbz_m )$ is generated by $e$ in dimension one
and $b$ in dimension two subject to
$e^2=kb$ and $b^{n+1} = eb^n = 0$. If $m$ is odd, we have to
change the first relation to $e^2= 0$.
\ele

\section{Surfaces}

\subsection{Orientable Surfaces}

In corollary \ref{Riemann} we have determined the homology of
$\sp{n}S$ for $S$ a Riemann surface of genus $g\geq 0$. This was based
on the construction of a chain complex for $\bsp{n}S$ based on cells
$e_1,\ldots, e_{2g}$ and $D$.  The coproduct at the chain level is
such that the $e_i$'s are primitive and
\begin{equation}\label{secondrelation}
D\mapsto D\tensor 1 + \sum e_i\tensor e_{i+g} -\sum e_{i+g}\tensor e_i
+ 1\tensor D\ \ \ \ \ \ 1\leq i\leq g
\end{equation}

If $f_1 = e_1^*,\ldots, f_{2g} = e_{2g}^*$ are the dual
cohomology classes, then their cup product satisfies
$f_if_{i+g} = b$, where $b$ as before is dual to the orientation class $D$.
There are no other relations in the cohomology of $S$.

Consider the diagram in lemma \ref{coproduct}. We propose to determine
for $i \neq j$
$\lambda_*(e_i\tensor e_j)\in C_*(X^2)^{\tensor 2}$. Here $\pi_*
(e_i\tensor e_j) = e_ie_j$. The effect $\chi^{-1}\delta_*^2$ on
$e_i\tensor e_j$ is as follows
\begin{eqnarray*}
e_i\tensor e_j&\mapsto&
\chi^{-1}\left(
(e_i\tensor 1\tensor 1\tensor 1 + 1\tensor e_i\tensor 1\tensor 1)
(1\tensor 1\tensor e_j\tensor 1 + 1\tensor 1\tensor 1\tensor e_j)\right)\\
&&=e_i\tensor e_j\tensor 1\tensor 1 + e_i\tensor 1\tensor 1\tensor e_j
- 1\tensor e_j\tensor e_i\tensor 1 + 1\tensor 1\tensor e_i\tensor e_j
\end{eqnarray*}
Applying $\pi_*$ to the left and
$\pi_*\tensor\pi_*$ to the right of this expression we obtain the coproduct
\begin{equation}\label{firstrelation}
e_ie_j\mapsto
e_ie_j\tensor 1 + e_i\tensor e_j - e_j\tensor e_i + 1\tensor e_ie_j
\end{equation}
Both of (\ref{firstrelation}) and (\ref{secondrelation}) pass to the
coproduct in homology. We will use throughout the same
symbol for a cycle in the chain complex and the homology class it generates.

By looking at dual classes in (\ref{firstrelation}) and
(\ref{secondrelation}) and pulling back, we see right away that
$e_i^*e^*_{i+g} = (e_ie_{i+g})^* + b$; i.e.

\ble\label{mainrelation} $(e_ie_{i+g})^* = e_i^*e_{i+g}^* - b$.
\ele

It turns out that this relation together with truncation by filtration
degree generate all relations in the cohomology of $\sp{n}S$.
The original calculation of $H^*(\sp{n}S)$ in
\cite{macdo} seems somewhat long winded. It can be
phrased in the following easier way. Note that since
$D$ generates a divided power algebra in $H_*(\spy (X))$
(corollary \ref{divided}), then its dual
$b$ generates a polynomial algebra $\bbz [b]$. On the other hand, the
dual of an exterior algebra is exterior and hence
\begin{equation}\label{infinitesym}
H^*(\spy (S))\cong E(f_1,\ldots, f_{2g})\tensor\bbz [b]
\end{equation}
Since $H_*(\sp{n}S)\lrar H_*(\spy S)$ is injective, then
$H^*(\spy S)\lrar H^*(\sp{n}S)$ is surjective, and
$H^*(\sp{n}S)$ is a quotient of (\ref{infinitesym}) by some
relations.

Consider the \emph{MacDonald relation} on the cohomology classes
$f_1,\ldots, f_{2g}$ and $b$
\begin{quote}
If $i_1, \dots ,i_a, j_1,\dots ,j_b,k_1,\dots, k_c$ are distinct
integers from 1 to $g$ inclusive, then provided that $a+b+2c+q = n+1$ we have
$$
f_{i_1}\cdots
f_{i_a}f_{j_1+g}\cdots f_{j_b+g}(f_{k_1}f_{k_1+g}- b)
\cdots (f_{k_c}f_{k_c+g} - b) b^q = 0$$
\end{quote}

\bth\cite{macdo}\label{macdocalc} $H^*(\sp{n}S)$ is the quotient of
$E(f_1,\ldots, f_{2g})\tensor \bbz [b]$ by the Macdonald relation.
\end{theorem}

\begin{proof}
We outline an alternative proof based on theorem \ref{main1} and lemma
\ref{coproduct}.  We know that $H_*(\sp{n}S)$ is rationally
 generated by $m$ fold
products of generators in $H_*(S)$, $m\leq n$. It follows that any
element $x$ of filtration degree $m\geq n+1$ cannot be in the image of
$\iota_* : H_*(\sp{n}X)\lrar H_*(\spy X)$ and hence $\iota^*(x^*) =
0\in H^*(\sp{n}S)$.  Choose a generator of filtration $(n+1)$ which we
write in the form
$$x= e_{i_1}\cdots e_{i_k}\sp{t}S\ \ \ , \ \ i_j\neq i_l, \ t+k = n+1$$
We can show by writing coproduct formulae that if no
pair of the form $\{e_s, e_{s+g}\}$ figures among the $e_{i_j}$'s
above, then the dual class verifies
$$(e_{i_1}\cdots e_{i_k}\sp{t}S)^* = f_{i_1}\cdots f_{i_k}b^t$$
where again the $f$'s are dual to the $e$'s.
If say $e_{i_1}=e_s$
and $e_{i_2}=e_{s+g}$, then
$(e_{i_1}\cdots e_{i_k}\sp{t}S)^* =
(e_se_{s+g})^*(e_{i_3}\cdots e_{i_k}\sp{t}S)^*$, and
$(e_se_{s+g})^*$ is as in lemma \ref{mainrelation}. In light of this,
the condition
$\iota^*\left(e_{i_1}\cdots e_{i_k}\sp{t}S\right)^* = 0$ translates
directly to the MacDonald's relation and this is the only such relation.
\end{proof}

\subsection{Non orientable Surfaces}
As far as we know the non-orientable case is not in the
literature. Let $U$ be the non-orientable surface of genus $g$. Then
$U$ is the connected sum of $g$-copies of the real projective plane
\begin{equation}\label{U}
U := \bbr\bbp^2\#\bbr\bbp^2\#\cdots\# \bbr\bbp^2\ \ \ \ \ \ \hbox{($g$-times)}
\end{equation}
We can write $U$ as a wedge of $g$-circles with a single
disc $D^2$ attached along the sum of
 degree two maps on each leaf.
If we denote as before the cellular generators by
$e_1,\ldots, e_g$ and $D\in C_2(U)$, then
\begin{equation}\label{boundary2}
\partial D = 2e_1 + \cdots + 2e_g
\end{equation}
The homology $H_*(\sp{n}U;\bbz )$ is completely determined by theorem
\ref{main1}. In particular of course $H_1(U) = \bbz_2\oplus\bbz^{g-1}$
and $H_2(U) = 0$.

\bre Note that $U$ has as oriented two cover a Riemann surface $S$ of
genus $g-1$.  This covering is obtained by embedding $S$ in $\bbr^3$
so that the origin $O$ is a center of symmetry.  The central symmetry
with respect to $O$ is a $\bbz_2$ free action on $S$ and the quotient
is $U$.  \ere

To determine the cohomology ring, we proceed as in the orientable case
and write the coproduct in $C_*(U)$.  Evidently the $e_i$'s are
primitive while
\begin{equation}\label{coproduct2}
D\mapsto D\tensor 1 + e_1\tensor e_1 + \cdots + e_g\tensor e_g+ 1\tensor D
\end{equation}
Since $D$ is a cycle modulo two, (\ref{coproduct2}) gives the
coproduct in homology modulo two as well. Moreover the coproduct on
classes $e_ie_j$ is as in (\ref{firstrelation}). Since
$H_2(\sp{n}U;\bbz_2)$ is generated by $D$ and $e_ie_j, i < j$, both
formulae yield the relation

\ble $b = f_1^2 = \cdots = f_g^2$.
\ele

We now show that this is the only relation in $H^*(\sp{n}U;\bbz_2)$
together with the filtration relation which demands that all
$n+1$-products be trivial.

\ble In $H^*(\sp{n}U;\bbz_2)$, 
$(e_{i_1}\ldots e_{i_r}\sp{t}D)^* = f_{i_1}\ldots f_{i_r}b^t$
where $r+t\leq n$ and $i_l\neq i_j$ for $l \neq j$.
\ele

\begin{proof}
Let's go back to the integral chain complex,
$\pi_* : C_*(U^n)\rightarrow C_*(\bsp{n}U)$.
We have $\pi_*(e_{i_1}\otimes \ldots \otimes e_{i_r}\otimes
D\otimes\cdots\otimes D) = t! e_{i_1}\ldots e_{i_r}\sp{t}D$.
Writing the coproduct for
this general class is notationally very involved. We get
precise enough of an idea by working out the case of
$e_1\sp{2}D$. We first write down $\delta_*^3 : C_*(X^3)\rightarrow
C_*(X^6)$ on $e_1\tensor D\tensor D$. This
consists of 18 terms obtained from the product
\begin{eqnarray*}
&& (e_1\tensor 1^{\tensor 5} + 1\tensor e_1\tensor 1^{\tensor 4})\\
&&(1^{\tensor 2}\tensor D\tensor 1^{\tensor 3} +
\sum 1^{\tensor 2}\tensor e_i\tensor e_i\tensor 1^{\tensor 2}
+ 1^{\tensor 3}\tensor D\tensor 1^{\tensor 2})\\
&&(1^{\tensor 4}\tensor D\tensor 1 +
\sum 1^{\tensor 4}\tensor e_i\tensor e_i
+ 1^{\tensor 5}\tensor D)
\end{eqnarray*}
As in the example of \S 3.1 we develop this expression,
shuffle by $\chi^{-1}$ and then
apply $\pi_*\tensor\pi_*$ to all terms to obtain the
coproduct in $C_*(\bsp{3} U)$
\begin{eqnarray*}
e_1\sp{2}D&\fract{\lambda_*}{\longmapsto}&
e_1\sp{2}D\tensor 1+ e_1\tensor\sp{2}D + \sp{2}D\tensor e_1
+1\tensor e_1\sp{2}D \\
&&+ \sum e_1e_iD\tensor e_i + e_1D\tensor D +
e_1e_i\tensor e_iD + \sum_{i<j} e_1e_ie_j\tensor e_ie_j\\
&& + \sum e_iD\tensor e_1e_i + D\tensor e_1D +
e_i\tensor e_1e_iD + \sum_{i<j} e_ie_j\tensor e_1e_ie_j 
\end{eqnarray*}
Reducing modulo two, and since all classes involved represent homology
classes, we arrive at the coproduct on $e_1\sp{2}D$
in $H_*(\sp{3}U;\bbz_2)$. The class $e_1\sp{2}D$ is the only basis element
whose image under $\lambda_*$ involves the basis element
 $e_1D\tensor D$. Consequently
\begin{equation}\label{exx}
(e_1\sp{2}D)^* = \lambda^*(e_1D\tensor D)^* =
 \lambda^*((e_1D)^*\tensor D^*) = \lambda^*(f_1b\tensor b):=f_1b^2
 \end{equation}
The general case is completely analogous.\end{proof}

Very much as in the proof of
theorem \ref{macdocalc} we can now deduce

\bpr\label{cohstructure} $H^*(\sp{n}U;\bbz_2)$ is generated by classes
$f_1,\ldots, f_g$ and $b$ under the relations\\
(i) $f_i^2=b$\\
(ii) $f_{i_1}\ldots f_{i_r}b^t = 0$,
$r+t= n+1$, $i_l\neq i_j$ for $l \neq j$.
\epr

\bex
Suppose $g=1$ and $U=\bbr\bbp^2$. Then 
$H^*(\sp{2}U;\bbz_2 )\cong\bbz_2[f_1]/(f_1^5)$
 which is the
cohomology of $\bbr\bbp^4$ (see lemma \ref{dupont}).
\eex

\bre\label{dt}
We can invoke the theorem of Dold and Thom to the effect that
$$\spy (X)\simeq \prod K(\tilde H_i(X;\bbz ),i)$$
for any finite type connected CW complex $X$.
Applying this to $S$ and $U$ we find
$$\spy (S)\simeq (S^1)^{2g}\times\bbc\bbp^{\infty}\ \ ,\ \
\spy (U)\simeq (S^1)^{g-1}\times\bbr\bbp^{\infty}$$

This is well in accordance with our homological calculations
since from proposition \ref{cohstructure} we can deduce that
$H^*(\spy U;\bbz_2)\cong E(h_1,\ldots, h_{g-1})\tensor\bbz_2 [c]$
where $h_i :=f_1+f_i$ and $c :=f_1$.
\ere


\section{Clifford's theorem and Analogs}

Let us consider 
the case when $S=C$ is a smooth complex projective curve,
or equivalently a compact Riemann surface. One then defines \cite{acgh}
the $n$-th \emph{Abel-Jacobi map} which is a holomorphic map
$$\mu_n : \sp{n}(C)\lrar J(C)$$
where $J(C)$, the ``Jacobian'' of $C$, is a complex torus of
dimension the genus of $C$. The maps $\mu_n$ are additive in the sense that
the following commutes
$$\begin{matrix}
\sp{r}(C)\times \sp{s}(C)&\lrar& \sp{r+s}(C)\cr
\decdnar{\mu_r\times\mu_s}&&\decdnar{\mu_{r+s}}\cr
J(C)\times J(C)&\lrar&J(C)\cr
\end{matrix}$$
the bottom map being addition in the abelian torus $J(C)$, and
the top map concatenation of points.
If $C$ is an elliptic curve for example,
then $J(C)\cong C$. The inverse preimages of $\mu_n$ are complex
projective spaces $\bbc\bbp^m$, where $m$ is related to the dimension of
some complete linear series on $C$(cf. \cite{acgh}). The dimension
of the preimages $\mu_n^{-1}(x), x\in J(C)$
is an upper semi-continuous function of $x$. The following is classical.

\bth\label{clifford} \cite{acgh}
 Write $\mu_n : \sp{n}C\rightarrow J(C)$,
$n\geq 1$, $g$ genus of $C$.\\
(i) If $n < 2g$, and $y \in \sp{n}C$, then
$\mu_n^{-1}(\mu_n(y)) = \bbc\bbp^{m(y)}$ for some $m(y) \leq {n\over 2}$. \\
(ii) If $n\geq 2g$ then $\mu_n^{-1}(x)= \bbc\bbp^{n-g}$ for all $x$.
\end{theorem}

The bounds in the theorem are sharp and are attained for hyperelliptic
curves.  Notice that part (ii) has a much more elaborate version due to
Mattuck and asserting that $\mu_n$ is a projectivized analytic bundle
projection with fiber $\bbc\bbp^{n-g}$.  The second part of this theorem
is due to Clifford and the proof is algebro-geometric in nature. We
now make the observation that theorem \ref{clifford} is in fact a
purely topological statement.

\bpr\label{clifford1} (Clifford's theorem : topological version)
Let $S$ be a closed oriented topological surface of genus $g>0$,
and $\bbc\bbp^m\lrar \sp{n}S$ a map that is non-trivial in homology.
\\ (i) If $n< 2g$, then necessarily $m\leq
\left[{n\over 2}\right]$;\\
(ii) if $n\geq 2g$, then $m\leq n-g$.
\epr

\begin{proof}
Let $h:\bbc\bbp^m\lrar\sp{n}S$ be a map such that
$h_*[\bbc\bbp^k]\neq 0$ for some $1\leq k\leq m$, where
$[\bbc\bbp^k]$ is the generator of $H_{2k}(\bbc\bbp^m )$. This says that
if $u\in H^2(\bbc\bbp^m)$ is the generator, then
there must be a class $x\in H^{2k}(SP^n(C))$ such that $h^*(x) = u^k$.
But by the cohomology structure of $SP^n(C)$ (theorem \ref{macdocalc}),
$x$ is decomposable into
a product of one dimensional generators and a single two dimensional
class $b$. Write $I$ the ideal generated by the one dimensional
generators. Since $h^*(x) = u^k\neq 0$, necessarily $x$ is
 $\pm b^k$ modulo terms in $I$, and hence $h^*(b) = \pm u$.
 Now $h^*(f_i)=0$ and
hence
$$\pm u^{[\frac{n}2] + 1} = h^*(b^{[\frac{n}2] + 1}) =
h^*( \prod_{k=1}^{[\frac{n}2] + 1} (b - f_{k}f_{k+g}) ) = 0$$
using the MacDonald relation (theorem \ref{macdocalc}).
This implies $m \leq [\frac{n}2]$.
One uses a similar argument for (ii).
\end{proof}

The next corollary recovers the original Clifford theorem.

\bco Choose a complex structure on $S$.
If $f: \bbc\bbp^m\lrar\sp{n}S$ is a non-constant holomorphic map, then
the conditions (i) and (ii) of proposition
 \ref{clifford1} hold.
\eco

\begin{proof}
It suffices to argue that a holomorphic map $f:
\bbc\bbp^m\lrar\sp{n}(S)$ is trivial in homology if and only if it is
constant.  Assume then that $f_*$ is trivial, and choose an embedding
$e:\sp{n}(S)\hookrightarrow\bbc\bbp^N$ realizing $\sp{n}S$ as a
projective variety \cite{acgh}. The composite $e\circ f:
\bbc\bbp^m\lrar\bbc\bbp^N$ is also trivial in homology. If we show
that $e\circ f$ is necessarily constant, then since $e$ is an
embedding it follows that $f$ is constant as well and hence the claim.

Let $g:\bbc\bbp^m\lrar\bbc\bbp^N$ be a non-constant holomorphic map.
Then $g$ is a finite ramified covering over its image $g(\bbc\bbp^m)$
which is a subvariety of $\bbc\bbp^N$. If $g(\bbc\bbp^m)$ is not
reduced to point, then it has dimension $\geq 1$ and its fundamental
cycle is non-trivial in $H_*(\bbc\bbp^N)$.  This fundamental cycle is
covered by a non-zero homology class in $\bbc\bbp^m$ (by a transfer
argument over $\bbq$ for example), and hence $g_*$ cannot be
trivial. So if $g_*$ is trivial, $g(\bbc\bbp^m)$ must be reduced to
point and $g$ is constant.
\end{proof}

Similarly there is a ``real analog" of proposition \ref{clifford1}
for \emph{unoriented} surfaces.

\bpr\label{clifford2} (Clifford's theorem : real version)
Let $U$ be a closed non-orientable topological surface of genus $g>0$,
and $\tau : \bbr\bbp^m\lrar \sp{n}U$ a map that is non-zero on
homology.  Then $m\leq 2n-g+1$.
\epr

\begin{proof}
We take the description in lemma \ref{cohstructure}. Write
$H^*(\bbr\bbp^m;\bbz_2) =\bbz_2 [u]/u^{m+1}$. Under the hypothesis,
and proceeding as in the first part of the proof of proposition
\ref{clifford1}, there is $f_i$ such that $\tau^*(f_i) = u$. Since
$f_1^2=\cdots =f_g^2=b$, we have also that $\tau^*(f_1)=\cdots =
\tau^*(f_g) = u$. Note that $f_1\cdots f_gb^{n-g+1} = 0$ in
$H^*(\sp{n}U;\bbz_2)$. But
$$\tau^*(f_1\cdots f_gb^{n-g+1} ) = u^gu^{2n-2g+2} = u^{2n-g+2} = 0$$
which means necessarily that $m + 1\leq 2n-g+2$ and hence the claim.
\end{proof}

We can see that the bound in proposition \ref{clifford2} is best
possible since in the genus 1 case (i.e. $U=\bbr\bbp^2$, covered by
$S^2$) we have the following result \cite{dupont}.

\ble\label{dupont} There is a homeomorphism
$\sp{n}(\bbr\bbp^2)\cong\bbr \bbp^{2n}$. \ele

\begin{proof} Write
$\bbr\bbp^2 = S^2/_{<T>}$ where $T$ is the antipodal involution acting
on $S^2$.  Note that $T$ extends to an action on $\sp{2n}(S^2)$ by
acting componentwise.  Now $T$ has no fixed points which implies that
the fixed point set of the action on $\sp{2n}(S^2)$ is
$\sp{n}(S^2/_{<T>})=\sp{n}(\bbr\bbp^2)$.  We need analyze this fixed
point set.

First of all if we write $S^2=\bbc\cup\{\infty\}$, then the action of
$T$ is $T(z) = -1/\bar z$. If on the other hand we identify $S^2$ with
$\bbc\bbp^1$ then in homogeneous coordinates we have $T([a:b]) =
[-\bar b:\bar a]$.

More generally identify $\sp{n}S^2$ with $\bbc\bbp^n$ as in lemma
\ref{sym2}.  That is identify first $\bbc\bbp^n$ with polynomials of
degree at most $n$, modulo scalar multiplication, by sending
$[a_0:\ldots : a_n]$ to $a_0+a_1z+\cdots +a_nz^n$. We can then check
that the map
$$\sp{n}S^2\lrar\bbc\bbp^n\ ,\ <z_1,\ldots, z_n>\mapsto (z+z_1)\cdots
(z+z_n)$$ 
is a homeomorphism (cf. \cite{hatcher}, chapter 4). Note that if $z_i$
coincides with $+\infty$, then the factor $``z+\infty"$ is omitted
from the product.  The action of $T$ on $\sp{2n}S^2$ translates to an
action on polynomials $(z+z_1)\ldots (z+z_{2n})\mapsto (z-1/\bar
z_1)\ldots (z-1/\bar z_{2n})$.  If $(a_0 = z_1\ldots z_{2n}, a_1 =
\sum z_1\cdots\hat z_i\cdots z_{2n}, \ldots, a_{2n-1} = z_1+\cdots +
z_{2n}, a_{2n}=1)$ are the coefficients of $p(z)= (z+z_1)\ldots
(z+z_{2n})$, modulo scalar, then $(1, (-1)\bar a_{2n-1}$, $\ldots,
(-1)^{i}\bar a_i,\ldots, (-1)^{2n}\bar a_0)$ are the coefficients of
$Tp(z)$, modulo scalars as well. After identification with
$\bbc\bbp^{2n}$, the antipodal action in homogeneous coordinates
becomes
$$T([a_0 : \cdots : a_{2n}])\longmapsto
[\bar a_{2n}:\cdots : (-1)^i\bar a_i : \cdots : \bar a_0]$$
The fixed point set $F\subset\bbc\bbp^{2n}$ of this action consists of all
$[a_0 :\ldots : a_{2n}]$ such that $\bar a_i= (-1)^ia_{2n-i}$ up to
usual scalar multiplication.
By splitting into real and imaginary parts, it is
easy to see that $F$ is a copy of $\bbr\bbp^{2n}$.
\end{proof}


\section{The Dold-Thom Homotopy Splitting}

Finally we point out how our previous constructions can be used to
give an elementary proof of the Dold-Thom splitting of $\spy X$ into a
product of Eilenberg-MacLane spaces (remark \ref{dt}).

For $X$ a two dimensional complex, theorem \ref{main1} shows that the
homology of $\sp{n}X$ only depends on the homology of $X$. So set
\begin{equation}\label{hom}
H_2(X) = \bbz^b\ \ \ ,\ \ \ H_1(X) =
\bbz_{n_1}\oplus\cdots\oplus\bbz_{n_r}\oplus \bbz^a
\end{equation}

We can assume that $C_*(X)$ is the chain complex
$\bbz^{r+b}\fract{\partial}{\lrar}\bbz^{r+a}\lrar 0$ sending basis
elements $D_1,\ldots, D_{r+b}$ to basis elements $e_1,\ldots, e_{r+a}$
according to
\begin{equation}\label{simplifiedcomplex}
\partial D_j = n_je_j\ ,\ 1\leq j\leq r,\ \ \ \
\partial D_{j} = 0,\ r < j\leq r+b
\end{equation}

We argue that the chain complex of $\bspy X$ in theorem \ref{main1} is
a tensor product of a certain number of copies of the chain complexes
for $\bspy S^1=S^1$, $\bspy S^2 = \bbc\bbp^{\infty}$ and $L_{n_j}=
S^{\infty}/\bbz_{n_j}$ the infinite Lens spaces.

It is evident indeed that the free generators in dimension one
generate a subchain complex of which homology is that of $(S^1)^a$,
while the free generators in dimension two contribute the homology of
$(\bbc\bbp^{\infty})^b$.  The subchain complex generated by
$\sp{i}(D_j)$ and $e_j \sp{i}(D_j)$ for fixed $j$ and $i \geq 0$ has
the homology of $L_{n_j}$.  If $X$ satisfies (\ref{hom}), then $\bspy
(X)$ is homologous to $Y = (S^1)^a\times (\bbc\bbp^{\infty})^b \times
L_{n_1}\ldots \times L_{n_r}$. But $Y$ is a generalized
Eilenberg-MacLane space and it classifies the cohomology of $X$.  So
there is a map $X\lrar Y$ which extends to $\spy X$ since $Y$ is an
abelian monoid;
$$\Psi : \spy (X)\lrar (S^1)^a\times (\bbc\bbp^{\infty})^b\times
L_{n_1}\times\cdots\times L_{n_r}$$
By the very construction of $C_*(\bspy X)$ and the fact that
$\Psi_*$ is a multiplicative map, we readily see that
$\Psi$ induces a homology isomorphism.
Since both spaces involved are
monoids again, they have abelian fundamental groups and are simple.
The map $\Psi$ is necessarily a homotopy equivalence. We deduce

\bpr Let $X$ be a two dimensional complex, and $H_*(X)$ as in
(\ref{hom}). Then there is a homotopy equivalence
$$\spy X\simeq (S^1)^a\times (\bbc\bbp^{\infty})^b\times
L_{n_1}\times\cdots\times L_{n_r}$$
\epr


\addcontentsline{toc}{section}{Bibliography}
\bibliography{biblio}
\bibliographystyle{plain[10pt]}

\end{document}